\documentclass[a4paper,12pt]{article}
\usepackage{xypic,amssymb,geometry,amsmath,enumerate,latexsym}
\newtheorem{example}{Example}[section]
\newtheorem{defn}[example]{Definition}
\newtheorem{prop}[example]{Proposition}
\newtheorem{thm}[example]{Theorem}

\newtheorem{cor}[example]{Corollary}

\newenvironment{pf}{{\bf Proof:}}{\hfill$\blacksquare$\mbox{}\vspace{.5cm}}
\def\epsilon{\varepsilon}
\def\cmod{\mathsf{Cmod}}
\def\Cmod{\mathsf{CMOD}}
\def\Aut{\mathsf{Aut}}
\def\AUT{\mathsf{AUT}}
\def\END{\mathsf{END}}
\def\End{\mathsf{End}}
\def\A{\mathcal A}
\def\C{\mathcal C}
\def\D{\mathcal D}
\def\B{\mathcal B}
\def\E{\mathcal E}
\def\<{\langle}
\def\>{\rangle}
\def\leq{\leqslant}
\begin{document}
\title{\Large \bf  HOMOTOPIES AND AUTOMORPHISMS \\ OF CROSSED MODULES OF GROUPOIDS
\thanks{MSC2000 Classification: 20L05,18D15,18D35,55D15,55U99 \newline
Keywords: crossed modules, groupoids, automorphisms, actors,
2-crossed modules}}

\small{ \author  {Ronald Brown  \\ University of Wales  \\ School
of Informatics, Mathematics Division\\ Bangor, Gwynedd \\ LL57
1UT, U.K.
\\ email: r.brown@bangor.ac.uk \\ \and \.{I}lhan \.{I}\c{c}en
\\ University of  \.{I}n\"{o}n\"{u} \\ Faculty of Science and Art
\\ Department of
Mathematics
\\ Malatya,  Turkey \\ email: iicen@inonu.edu.tr     }}
\maketitle

\def\subs{\subseteq}

\begin{abstract}
We give a detailed description of the structure of the actor
2-crossed module related to the automorphisms of a crossed module
of groupoids. This generalises work of Brown and Gilbert for the
case of crossed modules of groups, and part of this is needed for
work on 2-dimensional holonomy to be developed elsewhere.
\end{abstract}

\section*{Introduction}
This is the first of several planned papers whose aim is to
develop contexts for 2-dimensional versions of holonomy
\cite{Br-I}.

The notion of holonomy that we will use is that in Aof-Brown
\cite{Ao-Br}, which exposes ideas of Pradines justifying  \cite[
Th\'eor\`eme 1]{Pr}. This starts with a {\em locally Lie groupoid}
$(G,W)$ where $G$ is a groupoid and $W$ is a topological space
such that as a set $Ob(G) \subs W \subs G$ and certain other
conditions hold so that the groupoid operations on $G$ are `as
smooth as possible' on $W$. The locally Lie groupoid is {\it
extendible} if there is a topology on $G$ making $G$ a Lie
groupoid in which $W$ is an open subspace. In general $(G,W)$ is
not extendible but under extra conditions on the existence of
certain local continuous coadmissible sections, there is a Lie
groupoid $Hol(G,W)$ again containing $W$ as an open subspace and
with a morphism $Hol(G,W) \to G$ which is the identity on $W$.
Then $Hol(G,W)$ is the `smallest' overgroupoid of $G$ with this
property.

Here a {\em coadmissible section} is a section $s:X \to G$ of the
target map $\beta$ such that $\alpha s$ is a bijection on $X$. We
are investigating 2-dimensional versions of these ideas, and so we
need 2-dimensional versions of groupoids and of  coadmissible
sections.

There are various useful 2-dimensional notions of groupoid: double
groupoid, 2-groupoid, crossed module of groupoids. In this paper
we will consider homotopies in the context of crossed modules of
groupoids, which are equivalent to 2-groupoids, and to edge
symmetric double groupoids with the extra structure of connection.
In \S 5 we briefly summarise the application to 2-groupoids. It is
not clear at this stage whether a direct proof in that context
would or would not be simpler than that given here. In any case,
theory has been developed and applied for crossed modules of
groupoids and crossed complexes which  has not been developed for
the corresponding 2-groupoids and globular $\infty$-groupoids,
particularly the theory and applications of free objects (see for
example \cite{Br-99}).

The aim is to replace the notion of coadmissible section by {\em
coadmissible homotopy}. Our main results develop this relation
between coadmissible homotopies and automorphisms of crossed
modules of groupoids, and start by showing that a coadmissible
homotopy is of the form $f \simeq 1$ where $f$ is an automorphism.

The study of automorphisms of a crossed module ${\mathcal M}=(\mu
:M \to P)$ of a group $P$ was initiated by J.H.C. Whitehead in
\cite{Wh1}, where he showed the relation with derivations $s: P
\to M$. These derivations occurred in his work on homotopies of
his `homotopy systems' in \cite{Wh2}.  Work was continued by Lue
\cite{Lue} and Norrie \cite{Nor}, leading to the notion of the
{\em actor crossed square}
\begin{equation}\label{actor}
  \xymatrix{M \ar [d] _{\mu} \ar [r] & Der^*(P,M) \ar [d] \\
  P \ar [r] & \Aut({\mathcal M})   }
\end{equation}
which should be thought of as the direct generalisation  of the
{\em actor crossed module of a group} $M \to \Aut(M)$. See also
Breen \cite{Breen} for further uses of crossed squares.  Related
to this crossed square is a 2-{\em crossed module}
\begin{equation}\label{2-cr-m}
  M \to P \ltimes Der^*(P,M) \to \Aut({\mathcal M}).
\end{equation}
We shall later identify the group $P \ltimes Der^*(P,M)$ with the
group of invertible {\em free derivations} $P \to M$.

A different context for these results was given in Brown and
Gilbert \cite{Br-Gil} in terms of the monoidal closed category of
crossed modules of groupoids.  The technique of relating the
structure of {\em braided regular crossed modules} which arose in
this situation to simplicial groups allowed a construction of a
set of equivalent algebraic models of homotopy 3-types. The
extension from groups to groupoids is necessary for both these
contexts  since the object set of the internal hom is just the set
of morphisms of crossed modules, and the object set of the
automorphism structure has a group structure related to the
fundamental group of the 3-type.

However for  geometric applications it is natural to seek a
similar family of results, in particular the 2-crossed module
\eqref{2-cr-m}, for a crossed module of groupoids ${\cal C}=
(C,G,\delta)$, leading to the 2-crossed module
\begin{equation}  \label{equ:2-crs} M_2({\cal C})\stackrel{\zeta}{\longrightarrow} FDer^*({\cal
C})\stackrel{\Delta}{\longrightarrow} \Aut({\cal C})
\end{equation}with action of $\Aut(\C)$ on the other groups and
{\it Peiffer lifting} $$\<\;,\;\>: FDer^*({\cal C})\times
FDer^*({\cal C})\to M_2({\cal C}).$$ The part involving $\Delta$
will be the main focus of our work on 2-dimensional holonomy.

This has interesting special cases. If $\cal C$ is essentially a
groupoid $G$, then $ FDer^*({\cal C})$ can be identified with the
group $M(G)$ of coadmissible sections of $G$ and we have the
crossed module $$\Delta : M(G) \to \Aut(G)$$ which has occurred
sporadically in the literature (see for example \cite[\S 9.4,
Ex.5,6]{Br-88}, which are attributed to \cite{Ch}). Extensions by
a group $K$ of the type of this crossed module (c.f.
\cite{Br-Hig1}) are relevant to lax actions of the group $K$ on
the groupoid $G$ applied by Brylinski in \cite{Bry}.

From the point of view of the possibilities of `higher order
symmetry', it is interesting to note that a set corresponds to a
homotopy 0-type; the automorphisms of a set yield a group, which
corresponds to a homotopy 1-type; the automorphisms of a group
yield a crossed module, which corresponds to a homotopy 2-type;
and the automorphisms of a crossed module yield a 2-crossed
module, which corresponds to a homotopy 3-type.

One of our expository problems is that although the theory of
\cite{Br-Gil} deals with this general case of crossed modules of
groupoids, the interpretation for the groupoid case is not given.

We give below the full definitions and proofs of the algebraic
structure except for the verification of the laws for the
Peiffer lifting of the 2-crossed module, because we believe it
will help the reader to see explicitly the algebra that is
involved, and this will also make this work largely independent
of the papers Brown and Higgins \cite{Br-Hig2} and Brown and
Gilbert \cite{Br-Gil}. In particular, this makes our work
independent of the equivalence between crossed  complexes and
$\omega$-groupoids which is used by Brown and Higgins in
\cite{Br-Hig2}. We quote from  p.2 of this paper on the formulae
for the tensor product: ``Given formulae (3.1), (3.11) and
(3.14), it is possible, in principle, to verify all the above
facts within the category of crossed complexes, although the
computations, with their numerous special cases, would be long.
We prefer to prove these facts using the equivalent category
$\omega$-Grd of $\omega$-groupoids where the formulae are simpler
and have clearer geometric content''.

Thus we have in the above carried out a portion of this
verification. For the verification of the laws for the braided
part of the structure, we are however using facts from Brown and
Gilbert \cite{Br-Gil}.

In section 5 we mention briefly analogous results for the category
of 2-groupoids, whose equivalence with the category of crossed
modules of groupoids is the 2-truncated case of the main result
of \cite{Br-Hig4}. A later paper will develop this theory in its
own terms.

\section{Crossed modules of groupoids}

We recall the definition of crossed modules of groupoids. The
basic reference is  Brown-Higgins \cite{Br-Hig4}.

The source and target maps of a groupoid $G$ are written $\alpha,
\beta$ respectively. If  $G$ is {\em totally intransitive}, i.e.
if $\alpha=\beta$, then we usually use the notation $\beta$. The
composition in a groupoid $G$ of elements $a,b$ with $\beta a =
\alpha  b$ will be written additively, as $a + b$.  The main
reason for this is the convenience for dealing with combinations
of inverses and actions.

\begin{defn}{\em
Let $G, C$ be groupoids over the same object set and let $C$ be
totally intransitive.
Then an {\bf action}  of $G$ on $C$ is given by a partially defined
function
\[ \xymatrix {       C\times G \ar@{o->}[r]& C }        \]
written $(c, a)\mapsto c^a$, which satisfies

\begin{enumerate}[\rm (i)]
  \item $c^a$ is defined if and only if $\beta (c) = \alpha (a) $, and then
$\beta (c^a) = \beta (a)$,
  \item  $(c_1 + c_2)^a = {c_1}^a + {c_2}^a$,
  \item $c_1^{a+b} = (c_1^{a})^{b}$ and  $c_1 ^{e_x} = c_1$
\end{enumerate}
for all $c_1, c_2 \in C(x, x)$, $a\in G(x, y) $, $b\in G(y, z)$. }
\end{defn}

\begin{defn}{\rm
A {\bf crossed  module of groupoids} \cite{Br-Hig4} consists of a
morphism $\delta : C\rightarrow G$  of groupoids $C$ and $G$
which is the identity on the object sets such that $C$ is totally
intransitive, together with an action of $G$ on $C$ which
satisfies
\begin{enumerate}[\rm (i)]
  \item $\delta (c^a ) = -a + \delta c +a$
  \item  $c^{\delta c_1} = -{c_1}+ c + c_1 $
\end{enumerate}
for $c, c_1\in C(x, x)$, $a\in G(x, y)$.}
\end{defn}
It is a {\em pre-crossed module} if (ii) is not necessarily
satisfied. A (pre-)crossed module will be denoted by ${\cal C} =
(C, G, \delta)$. A {\bf crossed module of groups} is a crossed
module of groupoids as above in which $C$, $G$ are groups
\cite{Wh1}.

The followings are standard  examples of crossed modules of
groups and of groupoids:

\begin{enumerate}[\hspace{-1em} 1)]
  \item Let  $H$ be a normal subgroup of a  group $G$ with $i : H\to G$
the inclusion. The action of $G$ on the right of $H$ by
conjugation makes $(H, G, i)$     into  a crossed module. This
generalises to the case $H$ is a totally intransitive normal
subgroupoid of the groupoid $G$.
  \item  Suppose $G$ is a group and $M$ is a  right $G$-module;
 let $0: M\to G$ be the
constant map sending  $M$ to the identity element of $G$. Then
$(M, G, 0)$ is a crossed module.
  \item   Suppose given a   morphism
\[      \eta: M\to N   \]
of left $G$-modules  and form  the semi-direct product $G\ltimes
N$. This is a group which  acts on $M$ via the projection from
$G\ltimes N$ to $G$. We define a morphism
\[      \delta: M\to G\ltimes N  \]
by $\delta (m)  = (1, \eta (m))$ where $1$ denotes the identity
in $G$. Then $(M, G\ltimes N, \delta )$ is a crossed module.
\item   Labesse in \cite{La} defines  a {\it crossed group}. It
has been pointed out to us  by Larry Breen  at Coimbra in 1999
that a crossed group is exactly a crossed module  $(C, X\rtimes G,
\delta) $ where $G$ is a group acting on the set $X$, and
$X\rtimes G$ is the associated actor groupoid; thus the
simplicial construction from a crossed group described by Breen
in \cite{La} is exactly the nerve of the crossed module (see
\cite{Br-99} for an exposition).
\item Let $(Y,Z)$ be a pair of topological spaces and let $X$ be a
subset of $Z$. Then there is a crossed module of groupoids
$(C,G,\delta)$ in which $G$ is the fundamental groupoid
$\pi_1(Z,X)$ on the set $X$ of base points and for $ x \in X$,
$C(x)$ is the relative homotopy group $\pi_2(Y,X,x)$.
\end{enumerate}

 Let
${\cal C}$, ${\cal C'}$ be crossed  modules. A {\em morphism} $f
:{\cal C}\to {\cal C'} $  of crossed modules is a triple
$f=(f_0,f_1,f_2)$ such that $(f_0,f_1) $ form a morphism of
groupoids $G \to G'$ and $f_2$ is a family of morphisms $C \to
C'$ such that the following diagrams commute:
$$ \xymatrix{
           C \ar [d]_{\delta}\ar [r]^{f_2} &  C'\ar[d] ^{\delta'}     \\
                  G\ar [r]_{f_1}          &   G'}
 \qquad \xymatrix{  C\times G\ar [r]^{f_2\times f_1}  \ar @{o->}
[d]     &    C'\times G'   \ar @{o->} [d]  \\ C \ar [r]_{f_2 } &
C' } $$ So we can define a category $\cmod$ of crossed modules of
groupoids.


\section{Homotopies and derivations}

In this section, we combine the notion of coadmissible section,
which is fundamental to the  work of Ehresmann \cite{Eh1}, with
the derivations which occur in Whitehead's account \cite{Wh1} of
automorphism of crossed modules and which are later developed by
Lue, Norrie, Brown-Gilbert \cite{Br-Gil,Lue,Nor}.

Whitehead in \cite{Wh2} explored homotopies of morphisms of his
`homotopy systems' and this was put in the general context of
crossed complexes of groupoids by Brown and Higgins in
\cite{Br-Hig2}. So we are exploring the implications of the idea
that a natural generalisation of the notion of coadmissible
section for groupoids is that of coadmissible homotopy for crossed
modules.
\begin{defn}{\em
Let $g: \C \to \D$ be  a morphism of crossed modules
$\C=(C,G,\delta),\;
\D=(D,H, \delta')$ with object sets $X,Y$. A {\em homotopy} on $ g$ is a
pair $(s,g)$ where $s=(s_0,s_1)$ and  such that $s_0 \colon X\to
H,$ $s_1:G\to D$ and
\begin{align*}
\beta (s_0 x)          &= g_0x,  \\ \beta (s_1 a)        &= \beta
g_1(a),
\\ s_1(a + b ) &= s_1(a)^{g_1b}  + s_1(b),\\
\end{align*}  whenever $x \in X, a,b \in G $ and $a+b$ is defined. Such a function $s_1$ is
called a $g$-derivation. }
\end{defn}
\begin{prop}\label{pr1} Given a homotopy as above, the formulae  \begin{align*}f_0(x)  &= \alpha
s_0 (x), \\ f_1(a)  & = s_0(\alpha a) + g_1a +\delta s_1(a)
-s_0(\beta a),\\ f_2(c) &= (g_2c+s_1\delta c)^{-s_0\beta (c)}
\end{align*}
for all $x \in X, a \in G, c \in C$ define a morphism $f: \C \to
\D$ of crossed modules.
\end{prop}
The proof is given in the Appendix.

We therefore call $(s,g)$ a homotopy $f \simeq g$. This situation
can be displayed as:  $$ \spreaddiagramrows{2pc}
\spreaddiagramcolumns{2pc}
\def\objectstyle{\ssize} \def\labelstyle{\ssize}
\diagram
        C  \dto_{\delta} \rto^{f_2,g_2} & D \dto^{\delta'} \\
        G  \urto|{s_1}\dto<0.5ex> \dto<-0.5ex> \rto^{f_1,g_1}
             & H \dto<0.5ex> \dto<-0.5ex>              \\
        X  \urto|{s_0}          \rto_{f_0,g_0}      & Y
\enddiagram
$$

\begin{defn}{\rm For a homotopy $(s,1):f \simeq 1 : \C \to \C$, we
also call $s$ a {\em free derivation} on $\C$, and write
$\Delta(s)$ for $f$.
  }\end{defn} Let $FDer({\cal C})$ be the set of free
derivations of ${\cal C}$.

\begin{prop} \label{der}
Let ${\cal C} = (C, G, \delta)$ be a crossed module over groupoids.
Then a monoid structure on  $FDer({\cal C})$ is defined by the multiplication
\begin{align*}
(s*t)_{\epsilon}(z) = \begin{cases}
 (s_0g_0(z))+t_0(z) & \text{if }  \epsilon=0,  z\in X, \\
 t_1(z)+ (s_1g_1(z))^{t_0(\beta z)} & \text{if }\epsilon=1,  z \in G(x, y)
\end{cases}
\end{align*}
for $s, t\in FDer({\cal C})$ and $ g=\Delta (t)$. Further
$\Delta (s*t) = \Delta (s)\circ \Delta (t)$ and $ \Delta (1)=1$,
so that $\Delta $ is a monoid morphism
\[  FDer({\cal C})\to \End({\cal C})     .    \] \hfill$\blacksquare$
\end{prop}
The proof is given in the Appendix.

Let $FDer^*({\cal C})$ denote the group of invertible elements of
this monoid $ FDer({\cal C})$. An invertible free derivation is
also called a {\bf coadmissible homotopy}.

\begin{thm}\label{aut}
Let $s\in FDer({\cal C})$  and let $f=\Delta (s)$. Then the
following conditions are equivalent: {\em
\begin{enumerate}[(i)]
  \item  $s\in FDer^*({\cal C})$,
  \item   $f_1\in \Aut(G)$,
  \item   $f_2\in \Aut(C)$.
\end{enumerate}}\hfill$\blacksquare$
\end{thm}
The proof is given in the Appendix.

\begin{thm} \label{pre-crs}
There is an action of $\Aut({\cal C})$ on $s\in FDer^*({\cal C})$
given by
\begin{align*}
s^f (z) = \left \{ \begin{array}{ll} {f_1}^{-1} s_0 f_0(z), \ \ z
\in X \\ {f_2}^{-1} s_1 f_1(z),   \ \ z \in G.
\end{array} \right.
\end{align*}
 which makes
\[    \Delta \colon FDer^*({\cal C}) \rightarrow \Aut({\cal C})   \]
a pre-crossed module.
\hfill$\blacksquare$
\end{thm}
The proof is given in the Appendix.

The fact that $\Delta : FDer^*({\cal C})\to \Aut({\cal C})$ is a
pre-crossed module is a special case of results of Brown-Gilbert
\cite{Br-Gil}, which applies the monoidal closed structure of the
category of crossed complexes introduced by Brown-Higgins in
\cite{Br-Hig2} (see our section 4). The detailed description of
$\Delta :FDer^*({\cal C})\to \Aut({\cal C})$ is carried out in
Brown-Gilbert \cite[Proposition 3.3]{Br-Gil} only for the case
${\cal C}$ is a crossed module of groups.

We can also relate $FDer^*(\C)$ and $Der^*(\C)$, the group of all
invertible derivations of $\C$, where a {\em derivation} is simply
a pair $(s_0,s_1)$ in which $s_0$ is $x \mapsto 1_x$. The
multiplication of $Der^*(\C)$ is
\[  (s_1*t_1)(a)= t_1 a+s_1(a+\delta t_1a). \]
For a groupoid $G$ over $X$ let $M(G)$ denote  the group of
coadmissible sections of $G$ with the usual multiplication
$(s_0*t_0)(x) =s_0(\alpha(t_0x))+t_0(x)$.  Define an action of
$t_0\in M(G)$ on $s_1\in Der^*(\C)$ by
$$({s_1}^{t_0})(a)=s_1(t_0(\alpha a) +a-t_0(\beta a))^{t_0\beta
a}.$$ In view of the definition of the multiplication on
$FDer^*(\C)$ we have
\begin{thm}
There is a natural isomorphism
\[  FDer^*(\C)\to M(G)\ltimes Der^*(\C). \]
\end{thm}

The next theorem begins the analysis of the next level of the
actor structure.
\begin{thm}\label{2sec}
Let ${\cal C} = (C, G, \delta)$ be a crossed module over
groupoids with object set $X$. Let  $M_2(\C)$ be the group of
sections $s_2:X\to C$ of $\beta$ under pointwise addition. Then
there is a crossed module \begin{equation*}  \label{equ:2-crs1}
M_2({\cal C})\stackrel{\zeta}{\longrightarrow} FDer^*({\cal C}).
\end{equation*}where
$\zeta :M_2(\C)\to FDer^*(\C)$ is defined by $s_2\mapsto (\delta
s_2, s_1)$ where $s_1: a \mapsto  -s_2(\alpha a)^{a}+s_2(\beta
a)$, and the action  of $FDer^*(\C)$ on $M_2(\C)$ is given by
\[       (s_2)^{(t_0, t_1)}:x \mapsto  {s_2}(x)^{t_0(x)}.  \]
\end{thm}
The proof is given in the Appendix.


\section{Braided regular crossed modules \\ and 2-crossed modules}

In this section our object is to give the explicit relationship
between braided regular crossed modules and $2$-crossed modules.
This indicates a possible further context for development of work
on holonomy \cite{Br-I}.

The following material can be found in Brown and Gilbert
\cite{Br-Gil}.

Here we find it convenient to write a crossed module $\A$ as
$(A_2,A_1,\delta)$ instead of $(C,G)$ -- that is, we write
$A_2,A_1 $ for $C,G$ respectively and we also write $A_0$ for
$Ob(G)$. A {\em monoid bimorphism} $$b: \A,\A \to \A$$ consists of
a family of maps $$b_{ij} : A_i \times A_j \to A_{i+j}$$ for
$i,j=0,1,2; i+j \leq 2$ satisfying the following axioms, in which
$b_{11}(a,b)$ is written $\{a,b\}$.

\begin{enumerate}[ \bf 3.1]
\item $b_{00}$ gives $A_0$ the structure of monoid with identity
element written $e$ and multiplication written $(x,y)\mapsto xy$.
\item For $i=1,2$, $b_{0i}$ and $b_{i0}$ give actions of this monoid on the
left and right of the groupoid $A_i$ written $(x, c)\mapsto x.c$
and $(c, y)\mapsto c.y$ and which commute $x.(c.y)=(x.c).y$.
Further each action of an element of $A_0$ preserves the groupoid
structure, so that $x.(a+b)=x.a + x.b, (a+b).y = a.y+b.y$ and
$\alpha,\beta$ are equivariant with respect to these actions.
\item The actions are compatible with the actions of $A_1$ on
$A_2$ in the sense that if $c\in A_2(x)$, $a\in A_1(x, y),$ and $z
\in A_0$ then
\[        z.(c^a) = (z.c)^{z.a} \in A_2(z.x), \]
\[         (c^a).z = (c.z)^{a.z} \in A_2(x.z); \]
\item The boundary morphism $\delta : A_2 \to A_1$ is equivariant
with respect to these actions of $A_0$.
\item  $\{a, b\}\in A_2((\beta a)(\beta b)), \{O_e, b\}= O_{\beta
b}, \{a, O_e\} = O_{\beta a};$
\item  $\{a, b+b'\} = \{a, b\}^{\beta a.b'}+\{a, b'\};$
\item  $\{a+a', b\} = \{a',b\}+ \{a, b\}^{a'.\beta b};$
\item  $\delta \{a, b\}= -(\beta a.b)-a.\alpha b+\alpha
a.b+a.\beta b;$
\item  $\{a, \delta c'\}= -(\beta a.c')+(\alpha a.c')^{a.y} $ if
$c'\in A_2(y);$
\item  $ \{\delta c, b\} = -(c.\alpha b)^{x.b}+ c.\beta b$ if
$c\in A_2(p);$
\item   $x.\{a, b\} = \{x.a, b\},$
\ $\{a, b\}.x = \{a, b.x\},$  \ $\{a.x, b\} =\{a, x.b\},$
\end{enumerate}
for all $a, a', b, b'\in A_1, c, c'\in A_2$ and $x, y\in A_0$.

The crossed module $\A$ with this structure is called a {\em
semiregular braided crossed module}, and it is called {\em
regular} if the monoid $A_0$ is a group.

This structure is closely related to one given by Conduch\'{e} in
\cite{Con}. Recall from  \cite{Con}  that a {\bf 2-crossed module}
consists, in the first instance, of a complex of $P$-groups
\[   L \stackrel{\partial}{\longrightarrow} M\stackrel{\partial}{\longrightarrow} P \]
(so that $\partial \partial=0$) and $P$-equivariant homomorphisms,
where the group $P$ acts on itself by conjugation, such that
\[ L\stackrel{\partial}{\longrightarrow} M\]
is a crossed module, where $M$ acts on $L$ via $P$. We require
that    $(l^m)^p = (l^p)^{m^p}$ for all $l\in L, m\in M$, and
$p\in P$.  Further, there is a function $\< , \>: M\times M\to L,$
called a {\it Peiffer lifting}, which satisfies the following
axioms:

$P_1:$  $\partial\<m_0, m_1\> = {m_0}^{-1} {m_1}^{-1}m_0
{m_1}^{\partial m_0},$

$P_2:$   $\<\partial l, m\> = l^{-1} l^m ,$

 $P_3:$   $\<
m,\partial l\> = l^{-m} l^{\partial m} ,$

$P_4:$  $\<m_0, m_1m_2\> = \<m_0, m_2\>\<m_0,
m_1\>^{{m_2}^{\partial m_0}}$,

$P_5:$  $\<m_0 m_1, m_2\> = \<m_0, m_2\>^{m_1}\<m_1,
m_2^{{\partial m_0}}\>,$

$P_6:$ $\<m_0, m_1\>^p = \<{m_0}^p, {m_1}^p \>,$

for all $p \in P, m_0,m_1,m_2 \in M, l \in L$.

\setcounter{example}{11}
\begin{thm}\label{Braid}{\em (Brown and Gilbert \cite{Br-Gil})} The categories
of braided regular crossed modules and of \ $2$-crossed modules
are equivalent.
\end{thm}

\noindent {\bf Outline proof}

Let ${\A} = (A_2, A_1, \delta)$ be a regular crossed module.
Denote by $K$ the costar in $A_1$ at the vertex $e\in A_0,$ that
is, $K =\{a\in A_1 : \beta a = e\}.$ Then $K$ has a  group
operation given for any $a, b\in K$ by
\[         a b  =  b + (a.\alpha b).   \]
The source map $\alpha : K\to A_0$ is a homomorphism of groups and
is $A_0$-equivariant relative to the biaction of $A_0$ on $A_1$.
Note that the new composition extends the group structure on the
vertex group $A_1(e)$ so that $A_1(e)$ is  a subgroup of $K$: it
is plainly the kernel of $\alpha$. Further, $A_0$ acts diagonally
on $K$: for all $a\in K$ and $p\in A_0$ we set $a^p = p^{-1}.a.p$.
(There should be no confusion with the given action of $A_0$ on
$A_2$ which  we denote in a similar way.) Then the homomorphism
$\alpha :K\to A_0$ is $A_0$-equivariant relative to the diagonal
action on $K$ and the conjugation action of the group $A_0$ on
itself. Now $A_0$  also acts diagonally on the vertex group
$A_2(e)$ and so we have a complex of groups
\[     A_2(e)\stackrel{\delta}{\longrightarrow}
K\stackrel{\alpha}{\longrightarrow}A_0   \]
in which $\delta$ and $\alpha$ are $A_0$-equivariant. We know that
$\delta:A_2(e)\to A_1(e)$ is a crossed module: we claim that $K$
acts on $A_2(e)$, extending the action of $A_1(e)\subseteq K,$ so
that $\delta :A_2\to K$ is a crossed module.

We define an action $(c, a)\mapsto c!a$ by $c!a = (c.\alpha a)^a$
where $c\in A_2(e)$ and $a\in K.$ This is indeed  a group action
and $\delta$ is $K$-equivariant. Moreover, the actions of $A_2(e)$
on itself via $K$ and by conjugation coincide, for $\delta
:A_2(e)\to A_1(e)$ is a crossed module and so for all $c, c'\in
A_2(e),$
\[    c!\delta c' = (c.\alpha(\delta c'))^{\delta c'}= (c.e)^{\delta c'} = -c' + c + c'.  \]
Therefore the map $\delta :A_2(e)\to K$ is a crossed module.
Further the action of $A_0$ on $A_2(e)$ is compatible with that of
$K$.

The final structural component of a $2$-crossed module that we
need is the Peiffer lifting, which is provided by the braiding.
For suppose that ${\A}$ has a braiding $\{ , \} : A_1\times A_1
\to A_2$. Then the map $K\times K\to A_2(e)$ given by $(a,
b)\mapsto \{a^{-1}, b\}!a =\<a, b\>$ is a Peiffer lifting.
Therefore we have the $2$-crossed module
\[     A_2(e) \to K \to A_0.  \]
The verification of the axioms is obtained in \cite{Br-Gil} by
showing that this is the Moore complex of a simplicial group
$S({\A})$ whose Moore complex is of length 2 and applying the
results of \cite{Con}.

Now we show how a $2$-crossed module give rises to a braided
regular crossed module. So we begin with a $2$-crossed module
\[  L\stackrel{\partial}{\longrightarrow} G\stackrel{\partial}{\longrightarrow} P \]
and construct from it, in a functorial way, a regular, braided
crossed module ${\A}=(A_2, A_1, \delta )$.

The group of object of $A_0$ is just the group $P$. The underlying
set of elements of $A_1$ is $G\times P$ with source and target
maps $\alpha (g, p)=\partial(g)p$ and $\beta(g, p)= p.$ The
groupoid composition in $A_1$ is given by $(g_1, p_1)+(g_2, p_2) =
(g_1g_2, p_2)$ if $p_1 =\partial (g_2)p_2.$ The underlying set of
elements of $A_2$ is $L\times P$ with composition $(l_1, p)+(l_2,
p) = (l_1l_2, p).$ The boundary map $\delta: A_2\to A_1$ is given
by $\delta (l, p)=(\partial l, p)$ and the action of $A_1$ on
$A_2$ is given by $(l, p)^{(g, q)} = (l^g, q)$ if $p =\partial (g)
q.$ This does define a crossed module over $(A_1, A_0)$ and a
biaction of $A_0$ on ${\A}$ is obtained if we define
\begin{align*}
  p.(g, q) &= (g^{p^{-1}}, pq),  (g, q).p = (g, qp),  \\
 p.(l, q) &= (l^{p^{-1}}, pq),  (l, q).p = (l, qp),
\end{align*}
where  $(g, q), (l, q)\in A_2$ and $p\in A_0 = P$ and therefore
${\A}$ with this biaction is regular. The braiding on ${\A}$
is given by
\[  \{(g_1, p_1), (g_2, p_2)\} = (\<{g_1}^{-1}, {g_2}^{p_1}\>^{g_1}, p_1p_2) \]
where  $\<\; ,\; \>:G\times G\to L$ is the Peiffer lifting. \hfill
$\blacksquare$

In the next section we apply this result to automorphisms of
crossed modules of groupoids.


\section{ AUT({\cal C}) and 2-crossed modules}
In \cite{Br-Hig2} the category of crossed complexes is given the
structure of a monoidal closed category and this induces such a
structure on the category $\cmod$ of crossed modules of groupoids.
So there is a tensor product $-\otimes-$ and internal hom
$\Cmod(-,-) $ and  for all crossed modules $\A,\B,\C$ a natural
isomorphism
\[ \theta: \cmod(\A\otimes \B, \C)\to \cmod(\A, \Cmod(\B, \C)),  \]
which, together with the associativity of the tensor product,
implies the existence in $\cmod$  of a natural isomorphism
\[ \Theta: \Cmod(\A\otimes \B, \C)\to \Cmod(\A, \Cmod(\B, \C)).  \]
Further, the bijection
\[ \theta:\cmod(\Cmod(\A, \B)\otimes \A, \B)\to \cmod(\Cmod(\A, \B), \Cmod(\A, \B)) \]
shows that there is a unique morphism
\[  \epsilon_\A: \Cmod(\A, \B)\otimes \A\to \B \]
such that $\theta (\epsilon_\A)$ is the identity on $\Cmod(\A,
\B)$; $\epsilon_\A$ is called the {\it evaluation morphism}. Then
for all crossed modules $\A, \B, \C$  there is a morphism
\begin{align*} (\Cmod(\B, \C)\otimes \Cmod(\A, \B))\otimes \A&\to \Cmod(\B,
\C)\otimes (\Cmod(\A, \B)\otimes \A)\\ &\to \Cmod(\B, \C)\otimes
\B\to \C .
\end{align*}
 This corresponds under $\theta$ to a morphism
\[  \gamma_{\A\B\C} : \Cmod(\B, \C)\otimes \Cmod(\A, \B)\to \Cmod(\A, \C) \]
which is called {\it composition}.

We write $\END({\C })$ for $\Cmod({\C }, {\C }).$ The terminal
object in $\cmod$ is written $*$. There is a morphism $\eta_{{\C
}}: *\to \END({\C })$ corresponding to the isomorphism
$\lambda:*\otimes {\C }\to {\C }.$ The main result we need is the
following \cite{Kel}.
\begin{prop} \label{end}
The morphism $\eta_{\C } $ and the composition
\[ \mu = \gamma_{{\C }{\C }{\C }}:\END({\C })\otimes \END({\C })\to \END({\C }) \]
make $\END({\C })$ a monoid in $\cmod$ with respect to $\otimes$.
\end{prop}

\begin{cor}
The crossed module $\END(\C)$ may be given the structure of
braided semiregular crossed module.
\end{cor}
\begin{pf}
As shown in \cite{Br-Hig2,Br-Gil}, for any crossed module $\A$, a
morphism $\A\otimes \A \to \A$ corresponds exactly to a bimorphism
$\A,\A \to \A$, and a monoid structure is equivalent to a monoid
bimorphism.
\end{pf}

\subsection{$\AUT({\cal C})$ and 2-Crossed Modules}

Let $\A= \AUT({\cal C})$ be the full subcrossed  module of $\E=
\END(\C)$ on the object set $A_0= \Aut({\cal C})$ of automorphisms
of the crossed module ${\cal C} = (C, G, \delta)$. Thus $A_0$ is
the group of units of $E_0$. Now an element of $A_2$ is a section
over an automorphism of ${\cal C} = (C, G, \delta)$ and consists
of a pair $(s_2, f)$ where $s_2$ is a section and $f\in A_0$. An
element of $A_1$ is a homotopy over an automorphism of ${\cal C}
= (C, G, \delta)$ and consists of a triple $(s_0, s_1, f)$ where
$s_0$ is a section, $f\in A_0$, and $s_1$ is an $f$-derivation
$G\to C$ such that the endomorphism $f^0$ of ${\cal C} = (C, G,
\delta)$ which gives the source object of $(s_0, s_1, f)$ is
actually an automorphism. Clearly $f^0$ is an automorphism of
${\cal C} = (C, G, \delta)$ if and only if
\[       g(a) = s_0(x)+f(a)+\delta s_1(a) -s_0(y)   \]
\[      g(c) = (f(c)+s_1\delta (c))^{-s_0(\beta (c)} \]
\[       g(x) = \alpha s_0(x)  \]
for all $a\in G(x, y)$, $c\in C(x),$ $x\in X$ defines an
automorphism of ${\cal C} = (C, G, \delta)$.

Clearly $\A$ inherits from $\E$ the structure of a regular,
braided, crossed module. To determine the biaction of $A_0$ and
the braiding we have to understand the composition map $\gamma$
explicitly. A direct calculation leads to the following
non-trivial components for the bimorphism determining $\gamma:$
\begin{align*}
A_0\times A_0\to A_0 &:  (f_1, f_2)\mapsto f_1 f_2, \\ A_0\times
A_1\to A_1  &:   (f_1, (s_0, s_1, f))\mapsto (f_1s_0, f_1s_1,
f_1f),\\ A_1\times A_0\to A_1 &:  ((s_0, s_1, f), f_2)\mapsto
(s_0f_2, s_1f_2, ff_2),\\ A_0\times A_2\to A_2 &: (f_1, (s_2,
f))\mapsto (f_1(s_2), f_1f),\\ A_2\times A_0\to A_2 &: ((s_2, f),
f_2)\mapsto (s_2, ff_2),\\ A_1\times A_1\to A_2 &:  ((s_0, s_1,
f), (t_0, t_1, f'))\mapsto (s_1t_0, ff').
\end{align*}
These maps give  biactions of $A_0$ on $A_1,A_2$ and a braiding
$\{\, , \}: A_1\times A_1\to A_2.$ The monoid structure on $A_0$
is the usual composition of maps.

\begin{thm} \label{thm:2-crs}
The regular crossed module $\A =\AUT({\cal C})$ corresponds via
the equivalence  of Theorem \ref{Braid} to the 2-crossed module
\[  M_2(\C)\stackrel{\zeta}{\longrightarrow} FDer^*({\cal C})
\stackrel{\Delta}{\longrightarrow} Aut({\cal C}) . \]
\end{thm}
\begin{pf}
The costar in the groupoid $A_1$ at the identity automorphism $I$
of ${\cal C}$ may be identified as a set with $FDer^*({\cal C})$
and the group structure is given by $(s_0, s_1)*(t_0, t_1) =
(s_0*t_0, s_1*t_1)$  as in Proposition \ref{der} and Theorem
\ref{aut}. The vertex group $A_2(I)$ is identified with the group
$M_2(\C)$ with $\zeta (s_2) = (\delta s_2, s_1)$ as required in
Theorem \ref{2sec}. Note that $\Aut({\cal C})$ acts on
$FDer^*({\cal C})$ by
\[       (s_0, s_1)^f = (f^{-1} s_0 f, f^{-1} s_1 f)   \]
proved in Theorem \ref{pre-crs} and on $M_2(\C)$ by ${s_2}^f =
f^{-1}s_2f$. The action of $FDer^*({\cal C})$ on $M_2(\C)$ is
simply  $s_2^{(t_0, t_1)} = s_2^{t_0}$ and the Peiffer lifting is
given by
\begin{align*}
\<(s_0, s_1), (t_0,t_1)\> & =  \{(s_0, s_1)^{-1},
(t_0,t_1)\}!(s_0, s_1)\\
                        & =  (\{({s_0}^{-1}, ({s_1}^{-1})^{{s_0}^{-1}}), (t_0, t_1)\}.\Delta (s_0, s_1))^{(s_0, s_1)}\\
                        & =  ({s_1}^{-1})^{{s_0}^{-1}}*(t_0)^{s_0}\\
                        & =  {s_1}^{-1}({s_0}^{-1} *t_0* s_0).
\end{align*}
\end{pf}

\section{2-groupoids}

A commonly used 2-dimensional version of a groupoid is a
2-groupoid (see for example \cite{Kel-St2,Mo-Sv}). This is in fact
a 2-truncated case of the $\infty$-groupoids defined in 1981 in
\cite{Br-Hig5}, where it is shown that the categories of
$\infty$-groupoids and of crossed complexes are equivalent. In
particular, the categories of 2-groupoids and of crossed modules
of groupoids are equivalent. Hence the results of the previous
sections may be transferred to the category of 2-groupoids. This
gives the analogue of theorem \ref{thm:2-crs}:

\begin{thm} A 2-groupoid $\cal G$ determines a 2-crossed module of
the form
\[  M_2({\cal G})\stackrel{\zeta}{\longrightarrow} FDer^*({\cal G})\stackrel{\Delta}{\longrightarrow}
\Aut({\cal G}) . \]

\end{thm}

The precise description of these objects in 2-groupoid terms will
be left to another paper.

\section{Appendix}

\noindent
{\bf Proof of Proposition \ref{pr1}:}

We have to show that $f_1$ and $f_2$ are groupoid homomorphisms
and $f_2(c^a)= f_2(c)^{f_1(a)}$, for $c\in C(x), a\in G(x, y)$.
\begin{align*}
f_1(a+b) & =  s_0(x)+g_1(a+b)+\delta s_1(a+b)-s_0(z)\\
         & =  s_0(x)+g_1a+g_1b+\delta(s_1(a)^{g_1(b)}+s_1(b))-s_0(z)\\
         & =  s_0(x)+g_1a+g_1b-g_1b+\delta s_1 (a)+g_1b+\delta s_1(b)-s_0(z),\  \ \mbox{by definition of} \  \ \delta \\
         & =  s_0(x)+g_1a+\delta s_1(a)-s_0(y)+s_0(y)+g_1b+\delta s_1 (b)-s_0(z)\\
         & =  f_1(a)+f_1(b)\\
f_2(c+c')& =  (g_2(c+c')+s_1\delta (c+c'))^{-s_0(x)}\\
         & =  (g_2c+g_2c'+ s_1(\delta c+\delta c'))^{-s_0(x)}\\
         & =  (g_2c+g_2c'+ s_1(\delta c)^{g_1\delta c'}+ s_1(\delta c'))^{-s_0(x)}\\
         & =  (g_2c+g_2c'-g_2c'+s_1\delta c+g_2c'+s_1\delta c')^{-s_0(x)}\\
         & =  (g_2c+s_1\delta c+g_2c'+s_1\delta c')^{-s_0(x)} \\
         & =  f_2(c) + f_2(c').\\
\intertext{ Let $c\in C(x), a\in G(x, y)$. Then $\beta (c^a) =
\beta a$, $\beta c^a = y$. So} f_2(c^a) & =  (g_2c^a+ s_1\delta
(c^a))^{-s_0(\beta c^a)= -s_0(y)} \\
         & =  (g_2c^a + s_1(-a+\delta c+a ))^{-s_0(y)}\\
         & =  (g_2c^a + s_1(-a)^{g_1(\delta c+a)}+s_1(\delta c)^{g_1a
}+s_1(a))^{-s_0(y)}\\
&\qquad \qquad \mbox{\hfill (since  $-(s_1(a))^{-g_1a+g_1\delta c+g_1a}=
(s_1(-a))^{g_1(\delta c+a)}$),}\\
         & =  ( g_2c^a -(s_1a)^{-g_1a+g_1\delta c+g_1a}+(s_1\delta
         c)^{g_1a} + s_1a)^{-s_0(y)}\\
         & =  (g_2 c^a -s_1a^{(\delta g_2(c^{a}))}+(s_1\delta c)^{g_1a}+s_1a)^{-s_0(y)}\\
         & =  (-s_1(a)+ g_2(c)^{g_1a} +(s_1 \delta c)^{g_1a}+s_1(a))^{-s_0(y)}\\
         & =  (-s_1(a) +(g_2c+s_1\delta c)^{g_1a }+ s_1(a))^{-s_0(y)}\\
         & =  (-s_1(a)+(f_2(c)^{s_0(x)})^{g_1a} + s_1(a))^{-s_0(y)} \\
         & =  (f_2(c))^{s_0x+g_1a+\delta s_1a-s_0y}\\
         & =  f_2(c)^{f_1(a)}.
\end{align*}
So $f$ is an endomorphism of ${\cal C}$.
\hfill$\blacksquare$

\noindent {\bf Proof of Proposition \ref{der}:} It is clear that
$\beta (s*t)_0(x) = x$ and $\beta (s*t)_1(a) = \beta (a)$.

We have to show that $(s*t)_1$ is a derivation map. Let $a\in
G(x,y)$, $b\in G(y, z)$. Then
\begin{align*}
(s*t)_1(a+b) & =  t_1(a+b)+(s_1g_1(a+b))^{t_0(z)}\\ & =
t_1(a)^b+t_1(b)+(s_1(g_1(a)+g_1(b))^{t_0(z)}\\ & =
t_1(a)^b+t_1(b)+(s_1(g_1(a))^{g_1(b)}+ s_1(g_1(b))^{t_0(z)}\\ & =
t_1(a)^b+t_1(b)+(s_1(g_1(a))^{t_0(y)+b+\delta t_1(b)-t_0(z)}+
s_1(g_1(b))^{t_0(z)}\\ & =
t_1(a)^b+t_1(b)+s_1(g_1(a))^{t_0(y)+b+\delta t_1(b)}+
s_1(g_1(b))^{t_0(z)}\\ & =
t_1(a)^b+t_1(b)+(s_1(g_1(a)^{t_0(y)+b})^{\delta t_1(b)}+
s_1(g_1(b))^{t_0(z)}\\ & =
t_1(a)^b+t_1(b)-t_1(b)+(s_1(g_1(a))^{t_0(y)+b}+ t_1(b) +
s_1(g_1(b)^{t_0(z)}\\ & =  t_1(a)^b+(s_1(g_1(a))^{t_0(y)+b}+
t_1(b)+ s_1(g_1(b))^{t_0(z)}\\ & = (t_1(a)+s_1(g_1(a)^{t_0(y)})^b+
t_1(b)+ s_1(g_1(b)^{t_0(z)}\\ & = (s*t)_1(a)^b+ (s*t)_1(b).\\
\intertext{ For the associativity property, let $u, s, t\in FDer
({\cal C)}$ and let $f=\Delta (s), g=\Delta (t), h=\Delta (u)$.
Then} (u_0*(s*t)_0)(x) & =  (u_0(f_0g_0(x))+(t*s)_0(x)\\
                 & =   u_0(f_0g_0)(x)+(s_0g_0)(x)+t_0(x)\\
                 & =   u_0(f_0(g_0(x))+s_0g_0(x)+t_0(x)\\
                 & =   (u* s)_0(g_0(x))+t_0(x)\\
                 & =   ((u*s)_0*t_0)(x)\\
\intertext{and}
(u_1*(s*t)_1)(a) & =  (s*t)_1(a)+u_1(fg(a))^{(s*t)_0(x)} \\
                 & =  t_1a + (s_1ga)^{t_0(y)}+(u_1 fg a)^{(s*t)_0(x)} \\
                 & =  t_1a + (s*u)_1(ga)^{t_0(y)} \\
                 & =  ((u*s)_1*t_1)(a).\\
\intertext{Let  $s, t\in FDer({\cal C})$ be as above  and let
$a\in G(x, y)$. Then}
\Delta ({s*t})_1(a)& =  (s*t)_0(x) + a + \delta (s*t)_1(a) -(s*t)_0(y) \\
               & =  s_0g_0(x) + t_0(x) + a + \delta (t_1(a) + s_1g_1(a))^{t_0(y)}-(s_0g_0(y)+t_0(y))\\
               & =  s_0g_0(x) + t_0(x) + a + \delta (t_1(a) + \delta (s_1g_1(a)^{t_0(y)})-(s_0g_0(y)+t_0(y))\\
               & =  s_0g_0(x) + t_0(x) + a + \delta t_1(a) -t_0(y) + \delta s_1g_1(a)+t_0(y)-t_0(y)-s_0g_0(y)\\
               & =  s_0 g_0(x) +  \delta_t(a) + \delta s_1 g_1(a) - s_0g_0(y) \\
               & =  \Delta_s(\Delta_t)(a) \\
               & =   \Delta_s \circ \Delta_t(a). \\
\intertext{ Let $c\in C(x), a\in G(x, y)$.}
\Delta ({s*t})(c) & =  (c + (s*t)_1(\delta (c))^{-(s*t)_0(\beta c)}\\
                & =  (c +t_1(\delta (c))+s_1g_1(\delta (c))^{t_0(x)})^{-(s*t)_0(x)}\\
                & =  (c +t_1(\delta (c))^{-(s*t)_0(x)}+s_1\delta g_2(c))^{t_0(x)-(s*t)_0(x)} \\
                & =  (c +t_1(\delta (c))^{-t_0(x)-s_0g_0(x)}+s_1\delta g_2(c))^{-s_0g_0(x)} \\
                & =  ((c +t_1(\delta (c))^{-t_0(x)}+s_1\delta g_2(c))^{-s_0g_0(x)} \\
                & =  (\Delta_t(c)) + (s_1\delta g_2 (c))^{-s_0g_0(x)}, \ \ \mbox{since} \ \Delta_t(c)=g_2(c),\\
                & =  \Delta_s \circ \Delta_t (c)
\end{align*}
So $\Delta (s*t) = \Delta (s)\circ \Delta (t)$.

The identity of $Fder^*(\C)$ is  $c= (c_0, c_1)$ defined by $
c_0(x) =1_x  \ \ \mbox{and} \ \ c_1(a)= 1 $ for $x\in X$ and $a\in
G$. It is easy to see that $\Delta(c)= I$. \hfill$\blacksquare$

\noindent
{\bf Proof of Theorem \ref{aut}:}

That (i)$\Rightarrow$(ii), (i)$\Rightarrow$(iii) follows from the
fact that $\Delta$ is a morphism to $End({\cal C})$. We next
prove (ii)$\Rightarrow$(i). Suppose then $f_1\in Aut( G)$. We
define ${s^{-1}}= ({s_0}^{-1}, {s_1}^{-1}).$

Let ${s_0}^{-1}:X\to G, {s_1}^{-1}:G\to C$  by
\[  {s_0}^{-1}(x) = -s_0({f_0}^{-1}(x)) \ \  \mbox{and} \  \  {s_1}^{-1}(a) = -s_1({f_1}^{-1}(a))^{{s_0}^{-1}(y)}.\]
Since  $\beta {s_0}^{-1}(x) = x$ and $\alpha {s_0}^{-1}(x) =f_0(x)$,
${s_0}^{-1}$ is
an inverse element of $s_0$. In fact,
\begin{align*}
           (s^{-1}* s)_0(x) & =  {s_0}^{-1}({f_0}^{-1})(x) + s_0(x) \\
                          &  =  -s_0({f_0}^{-1}({f_0})(x)) + s_0(x) \\
                          & =  -s_0(x) + s_0(x) \\
                          & =  c_0(x)=1_x
\end{align*}
and also
\begin{align*}
            ( s* s^{-1})_0(x) & =  s_0 ({f_0}^{-1})(x) + {s_0}^{-1}(x)\\
                            & =  s_0({f_0}^{-1})(x) - s_0({f_0}^{-1}(x) \\
                            & =  c_0(y)=1_y.
\end{align*}

We have to show that ${s_1}^{-1}$ is a derivation map. Let  $a, b,
a+b\in G$ and let $a'= {f_1}^{-1}(a), b'={f_1}^{-1}(b)$,  $\beta b
= z$. Note that  ${s_0}^{-1}{\beta(a+b)={s_0}^{-1}(z)}$ and
${-s_0(z)= -s_0{f_0}^{-1}(z)}.$
\begin{align*}
{s_1}^{-1}(a+b) & =  -(s_1{f_1}^{-1}(a+b))^{{s_0}^{-1}{\beta(a+b)}}, \ \ \mbox{by definition of}\  {s_1}^{-1}\\
                & =  -(s_1({f_1}^{-1}a+{f_1}^{-1}b))^{{s_0}^{-1}(z)} \\
                & =  -(s_1(a'+b'))^{{s_0}^{-1}(z)}\\
                & =  -((s_1a')^{b'}+ s_1(b'))^{{s_0}^{-1}(z)}, \ \ \mbox{since} \ s_1 \ \mbox{is a derivation},\\
                & =  (-s_1(b')-(s_1(a'))^{b'})^{-s_0(z)}\\
                & =  -(s_1(b')+(s_1(a'))^{b'+\delta s_1(b')})^{-s_0(z)}\\
                & =  -(s_1(b')^{-s_0(z)}+(s_1(a'))^{b'+\delta s_1(b')-s_0{f_0}^{-1}(z)}\\
\intertext{ Since $f(b') = s_0{f_0}^{-1}(y)+b'+\delta s_1(b')
-s_0{f_0}^{-1}(z)= b$ and $b'+\delta s_1(b') -s_0{f_0}^{-1}(z)=
b- s_0{f_0}^{-1}(y),$}
                & =  -(s_1(b')^{-s_0(z)}+(s_1(a'))^{b-s_0{f_0}^{-1}(y)}\\
                & =  -(s_1(b')^{-s_0(z)}+(s_1(a'))^{b-{s_0}^{-1}(y)}\\
                & = - (s_1(a'))^b)^{-{s_0}^{-1}(y)} -s_1(b')^{{s_0}^{-1}}\\
                & = - (s_1({f_1}^{-1}(a))^b)^{-{s_0}^{-1}(y)} -s_1{f_1}^{-1}(b')^{{s_0}^{-1}}\\
                & =  {s_1}^{-1}(a)^b+{s_1}^{-1}(b). \\
\intertext{ One can easily show that $s*s^{-1} = c$ and $s^{-1} *
s = c$.}
          (s*s^{-1})_1(a) & =  {s_1}^{-1}(a) + s_1 (f^{-1}(a))^{{s_0}^{-1}(y)} \\
                        & =  -s_1f^{-1}(a)^{{s_0}^{-1}(y)}+s_1 (f^{-1}(a))^{{s_0}^{-1}(y)}  \\
                        & =  c_1(a).
\end{align*}
and
\begin{align*}
        (s^{-1}*s)_1(a) & =   s_1(a) + {s_1}^{-1}(f(a))^{s_0(y)}\\
                       & =   s_1(a) - s_1(f^{-1}(f(a))^{s_0(y)})^{{s_0}^{-1}(y)} \\
                       & =  s_1(a)-s_1(a)  \\
                       & =  c_1(a).
\end{align*}
Now we will prove (iii)$\Rightarrow$(i). We first recalculate
$(s*t)_1$ in terms of $f_2$. Let $\Delta(t) = g$ and let
$\Delta(s) = f$, $a\in G(x, y)$ as above.
\begin{align*}
(s*t)_1(a) & =  t_1(a)+s_1 g_1(a)^{t_0(y)}\\
        & =  t_1(a)+s_1(t_0(x)+a+\delta t_1(a)-t_0(y))^{t_0(y)}\\
        & = t_1(a)+s_1(t_0(x))^{a+\delta t_1(a)-t_0(y)}+s_1(a)^{\delta t_1(a)-t_0(y)}
             +s_1(\delta t_1(a))^{t_0(y)})+{s_1(-t_0(y))}^{t_0(y)}\\
        & = t_1(a)+s_1(t_0(x))^{a+\delta t_1(a)}+s_1(a)^{\delta t_1(a)}
             +s_1(\delta t_1(a))+s_1(-t_0(y))^{t_0(y)}\\
        & = t_1(a)+s_1(t_0(x)^a)^{\delta t_1(a)}-t_1(a)+s_1(a)+t_1(a)
             +s_1(\delta t_1(a))-s_1(t_0(y)),
\end{align*}
since $s_1(-t_0(y))^{t_0(y)} = -s_1t_0(y)$,
\begin{align*}
        & = t_1(a)-t_1(a)+s_1(t_0(x)^a+t_1(a)-t_1(a)+s_1(a)+t_1(a)
             +s_1(\delta t_1(a))-s_1(t_0(y))\\
        & = s_1(t_0(x))^a+s_1(a)+t_1(a)
             +s_1(\delta t_1(a))-s_1(t_0(y))\\
        & = s_1(t_0(x))^a+s_1(a)+f_2(t_1(a))^{s_0(y)}-s_1(t_0(y))\\
\end{align*}
Now, suppose that $f_2$ has
inverse  ${f_2}^{-1}$. Let ${s^{-1}}= ({s_0}^{-1}, {s_1}^{-1})$ be defined by
\[    {s_0}^{-1}(x) = -s_0{f_0}^{-1}(x), \ \ x\in X \]
\[    {s_1}^{-1}(a) = {f_2}^{-1}(-s_1(a)-(s_1{s_0}^{-1}(x))^a+(s_1{s_0}^{-1}(y))^{-s_0(y)}, \ \ a\in G(x, y)  \]
We prove that $s^{-1}$ is an inverse element of $s$ and is a derivation map.
Clearly
\[  s*s^{-1}(x) = c_0(x) \ \  \mbox{and} \ \  s^{-1}*s(x) = c_0(x)   \]
by the argument as above.

Next we prove $(s*s^{-1})_1(a) = c_1(a)$,  for $a\in G(x, y)$.
\begin{align*}
(s*s^{-1})_1(a)  & =  (s_1({s_0}^{-1}(x))^a+s_1(a)+f_2({f_2}^{-1}((-s_1(a)+(s_1s_0{f_0}^{-1}(x))^a\\
                 &  \quad -(s_1s_0{f_0}^{-1}(y))^{-s_0(y)})^{s_0(y)})-s_1{s_0}^{-1}(y) \\
                 & =(s_1({s_0}{f_0}^{-1}(x))^a+s_1(a)-s_1(a)+(s_1s_0{f_0}^{-1}(x))^a  \\
                 &  \quad  +(s_1s_0{f_0}^{-1}(y))-s_1{s_0}{f_0}^{-1}(y) \\
                 & = c_1(a).
\end{align*}
Since $(s*s^{-1}) =1$ and also $s^{-1} *s' =1$.
It follows that
${s_1}^{-1}*{s_1} = {s_1}^{-1}*{s_1}*{s_1}^{-1}*{s_1}' = {s_1}'* ({s_1}*{s_1}^{-1})* {s_1}'= {s_1}^{-1}*{s_1}'= 1$
and so ${s_1}^{-1} = {s_1}'$, i.e.,
\[  (s^{-1}*s)_1(a)  = c_1(a),  \ \ \mbox{for all} \ \ \ a\in G.  \]

We have to prove that ${s_1}^{-1}$ is a derivation map. Let
$a\in G(x, y), b\in G(y, z)$. We write
$f_2{s_1}^{-1}(a) = (-s_1(a)-(s_1{s_0}^{-1}(x))^a+(s_1{s_0}^{-1}(y))^{-s_0(y)}$, and
then
\begin{align*}
f_2(({s_1}^{-1}(a))^b+{s_1}^{-1}(b))& =  f_2 ({s_1}^{-1}(a))^b+ f_2({s_1}^{-1}(b))\\
                                  & = (-s_1(b)+(f_2({s_1}^{-1}(a))^{s_0(y)})^b+s_1(b))^{-s_0(z)}+f_2({s_1}^{-1}(b)),\\
                                 & = (-s_1(b)^{-s_0(z)}+(-s_1(a)-(s_1{s_0}^{-1}(x))^a+s_1{s_0}^{-1}(y))^{b-s_0(z)}  \\
                                  &  \quad \ \ +s_1(b)^{-s_0(z)} (-s_1b)^{-s_0z}-(s_1{s_0}^{-1}(y)))^{b-s_0(z)}+s_1{s_0}^{-1}(z)^{-s_0(z)} \\
                                 & = (-s_1(b)^{-s_0(z)}-s_1(a)^{b-s_0(z)}-(s_1{s_0}^{-1}(x))^{a+b-s_0(z)}+s_1{s_0}^{-1}(z))^{-s_0(z)}     \\
                                 & = (-s_1(b)-s_1(a)^b-(s_1{s_0}^{-1}(x))^{a+b}+s_1{s_0}^{-1}(z))^{-s_0(z)}     \\
                                 & = (-(s_1(a)^b+s_1(b))-(s_1{s_0}^{-1}(x))^{a+b}+s_1{s_0}^{-1}(z))^{-s_0(z)}     \\
                                 & =  f_2({s_1}^{-1}(a+b))
\end{align*}
Hence ${s_1}^{-1}$ is a free derivation. \hfill$\blacksquare$

\noindent
{\bf Proof of Theorem \ref{pre-crs}:}

Now, we will show that
\[  FDer^*({\cal C})\times Aut({\cal C}) \rightarrow FDer^*({\cal C})    \]
\[     (s, f)\mapsto s^f              \]
is an action of $Aut({\cal C})$ on $FDer^*({\cal C})$.

\begin{align*}
s^f (z) = \left \{ \begin{array}{ll} {f_1}^{-1} s_0 f_0(z), \ \
z\in X \\ {f_2}^{-1} s_1 f_1(z),   \ \ z\in G.
\end{array} \right.
\end{align*}

In fact this give rise to an action over a groupoid:

\begin{align*}
s^{fg} (z) = \left \{ \begin{array}{ll} {(fg)_0}^{-1} s_0
(fg)_0(z) = {g_0}^{-1}{f_0}^{-1} s_0 f_0g_0 (z) =
{g_0}^{-1}{s_0}^{f_0}{g_0}(z) = ({s_0}^{f_0})^{g_0}(z), \  z \in X
\\ {(fg)_1}^{-1} s_1 (fg)_1(z) = {g_1}^{-1}{f_1}^{-1} s_1 f_1g_1
(z) = {g_1}^{-1}{s_0}^{f_1} g_1(z) = ({s_1}^{f_1})^{g_1}(z),   \
z\in G \\
\end{array} \right.
\end{align*}
and
\begin{align*}
s^I (z) = \left \{ \begin{array}{ll} I^{-1} s_0 I(z) = s_0(z),
z\in X \\ I^{-1} s_1 I(z) = s_1(z) ,  z\in G.
\end{array} \right.
\end{align*}
Let $s: f\simeq I$. Then $s^f : \bar{f} \simeq I$.
We can show $s^f$ as  the following  diagram:
$$
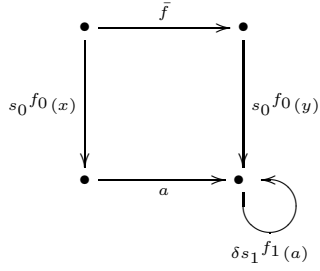
\spreaddiagramrows{2pc} \spreaddiagramcolumns{2pc}
\def \objectstyle {\ssize} \def \labelstyle{\sssize}
\diagram
\bullet \dto_{{s_0}^{f_0}(x)} \rto^{\bar{f}} & \bullet \dto^{{s_0}^{f_0}(y)} \\
\bullet \rto_a & \bullet \ \todr_{\delta {s_1}^{f_1}(a)}
\enddiagram
$$
Is $s^f\in FDer^*({\cal C})$ ?
Clearly one can see $\beta {s_0}^{f_0}{x} = x$ and $\beta {s_1}^{f_1}(a) = \beta (a).$
Also we should have to show that $s^f(a+b) = s^f(a)^b + s^f(b)$. We have
\[  \bar{f}(a+b)= {s_0}^{g_0}(x) + (a+b)+\delta s^f(a+b)-{s_0}^{f_0}(z) \]
by definition of $\bar{f}$ and
\begin{align*}
\bar{f}(a)+\bar{f}(b) &= {s_0}^f(x) +a+\delta {s_1}^f(a)+b+\delta {s_1}^f(b)-{s_1}^f(z) \\
                  &= {s_0}^f(x) +a+ b-b+\delta {s_1}^f(a)+b+\delta {s_1}^f(b)-{s_1}^f(z) \\
                  &= {s_0}^f(x) +a+b+\delta ({s_1}^f(a)^b + {s_1}^f(b) -{s_0}^f(z)  \\
                  &= \bar{f}(a+b)
\end{align*}
So  $s^f(a)^b +s^f(b) = s^f(a+b)$
and also we can obtain
\[ I(a) = -{s_0}(x) + \bar{f}(a) + {s_0}^f(y)- \delta s^f (a). \]
\begin{align*}
f^{-1} \Delta s f(a)&=    f^{-1}( {s_0}f(x) +f(a)+ \delta_s(fa) - {s_0}^f(y))\\
                    &=    f^{-1}{s_0}f(x) + f^{-1}f(a)+  f^{-1}\delta_s(fa) - f^{-1} {s_0}^f(y)) \\
                    &=  {s_0}^f(x) +a + \delta f^{-1}s_1 f(a) - {s_0}^f(y) \\
                    &=  {s_0}^f(x) +a + \delta {s_1}^f(a) - {s_0}^f(y) \\
                    &=  \Delta (s^f) (a).
\end{align*}
Hence $\Delta (s^f)(a) =  f^{-1} \Delta s f(a)$.
\hfill$\blacksquare$

\noindent
{\bf Proof of Theorem \ref{2sec}:}
 The group structure on $M_2(\C)$ is pointwise multiplication. If
 $\zeta(s_2)= (s_0,s_1)$ then $s_1$  is a derivation since
\begin{align*}
      s_1(a+b) & =  (-s_2(x))^{(a+b)} +s_2(z)\\
                       & =  (-s_2(x))^{(a+b)}+s_2(y)^{(b)}-s_2(y)^{(b)} +s_2(z)\\
                       & =  ((-s_2(x))^{(a)}+s_2(y))^{(b)}-s_2(y)^{(b)} +s_2(z)\\
                       & =  (s_1(a))^{(b)}+ s_1(b).
\end{align*}

The action of $FDer^*(\C)$ on $M_2(\C)$ is
\[       (s_2)^{(t_0, t_1)} = {s_2}^{t_0}:x \mapsto (s_2x)^{t_0x}  \]
and we have to prove
\begin{enumerate}[(i)]
  \item $\zeta((s_2)^{(t_0, t_1)})= (t_0, t_1)^{-1}*\zeta(s_2)*(t_0, t_1),$
  \item $(s_2)^{\zeta(t_2)} = (t_2)^{-1}*(s_2)*(t_2).$
\end{enumerate}
Then  $\zeta((s_2)^{(t_0, t_1)}= \zeta ({s_2}^{t_0})
=(\delta({s_2}^{t_0}),  t')$ say where if $a:x \to y$
\begin{align*}
                t'(a) & =  (-{s_2}^{t_0} (x))^{a}+{s_2}^{t_0}(y), \ \ a\in G(x, y)\\
                       & =  (-{s_2(x)}^{t_0(x)} )^{a}+{s_2(y)}^{t_0(y)},\\
                       & =  -{s_2(x)}^{{t_0} (x)+{a}}+({s_2}^{t_0})(y).
\end{align*}
On the other hand,
\begin{align*}
(t_0, t_1)^{-1}*\zeta(s_2)*(t_0, t_1)& =
({t_0}^{-1},{t_1}^{-1})*(\delta(s_2),s_1)*(t_0, t_1), \\
                   & =   ({t_0}^{-1}*\delta({s_2})*t_0, {t_1}^{-1}*s_1* t_1),\\
                   & =   (\delta({s_2}^{t_0}),  {t_1}^{-1}*s_1* t_1).
\end{align*}
So we have to show that $ ({t_1}^{-1}*s_1* t_1) (a)=
-{s_2(x)}^{{t_0} (x)+{a}}+({s_2}^{t_0})(y).$ Clearly $\zeta(s_2)*
(t_0, t_1) = (\delta(s_2), s_1)*(t_0, t_1)= (\delta(s_2)*t_0,
s_1*t_1).$ and $ (s_1* t_1)(a) = t_1(a)+s_1(a)=
t_1(a)+(-s_2(x)^{a}+s_2(y))^{t_0(y)}.$ Then
\begin{align*}
({t_1}^{-1}*s_1* t_1)(a) & =
(s_1*t_1)(a)+({t_1}^{-1}(a))^{(\delta(s_2)*t_0)(y)}\\
                          & = t_1(a)+(-s_2(x)^{a}+s_2(y))^{t_0(y)}+({t_1}^{-1}(a))^{(\delta(s_2)*t_0)(y)} \\
                          & = t_1(a)-s_2(x)^{a+t_0(y)}+s_2(y)^{t_0(y)}+(-{t_1}(a)^{-t_0(y)})^{(\delta(s_2)*t_0)(y)} \\
& =
t_1(a)-s_2(x)^{a+t_0(y)}+s_2(y)^{t_0(y)}+(-{t_1}(a))^{-t_0(y)+(\delta(s_2))*t_0(y)}
\\ & =
t_1(a)-s_2(x)^{a+t_0(y)}+s_2(y)^{t_0(y)}+(-{t_1}(a))^{(\delta({s_2}^{t_0})(y)}
\\ & =  t_1(a)-s_2(x)^{a+t_0(y)}+s_2(y)^{t_0(y)}-
s_2(y)^{t_0(y)}-t_1(a) +{s_2(y)}^{s_0(y)}\\ & =
-s_2(x)^{a+t_0(y)-\delta(t_1(a))}+ {s_2(y)}^{t_0(y)}\\ & =
-s_2(x)^{t_0(x)+a}+ {s_2}^{t_0}(y)
 \\ & =  t'(a).
\end{align*}
This proves (i).

To prove (ii) we note that
\begin{align*}
(s_2)^{\zeta (t_2)} & =  (s_2)^{(\delta(t_2), t_1)}\\
                          & =  ({s_2}^{\delta(t_2)})\\
                          & =  ({t_2}^{-1}*s_2*t_2)
\end{align*}
as is required.

{}

\end{document}